\def\Z{{Z\!\!\!\! Z}}
\def\CR{\hbox{{$\cal R$}}}\def\CC{\hbox{{$\cal C$}}}

\def\tens{\mathop{\otimes}}

\def\id{{\rm id}}

\def\proof{{\bf Proof\ }}
\def\endproof{{\ $\diamond$}\medskip}
\def\eqn#1#2{\begin{equation}#2\label{#1}\end{equation}}

\documentstyle{article}
\newtheorem{propos}{{\rm PROPOSITION}}
\newtheorem{corol}[propos]{{\rm COROLLARY}}

\newtheorem{lemma}[propos]{{\rm LEMMA}}
\newtheorem{example}[propos]{{\rm EXAMPLE}}

\begin{document}

{\ }\qquad  \vspace{.2in}

\parindent 0pt

{\Large $\Z_n$-QUASIALGEBRAS}

\baselineskip 13pt
\bigskip
\bigskip
HELENA ALBUQUERQUE\footnote{Supported by CMUC-JNICT and by Praxis
2/2.1/Mat7458/94}

\medskip
Departamento de Matematica-Faculdade de Ciencias e Tecnologia

Universidade de Coimbra, Apartado 3008

3000 Coimbra, Portugal

\bigskip
SHAHN MAJID\footnote{Royal Society University Research Fellow and
Fellow of Pembroke College, Cambridge}

\medskip
Department of Applied Mathematics and Theoretical Physics

University of Cambridge

Cambridge CB3 9EW, UK

\bigskip
\medskip ABSTRACT Recently we have reformulated the octonions
as quasissociative algebras (quasialgebras) living in a symmetric
monoidal category. In this note we provide further examples of
quasialgebras, namely ones where the nonassociativity is induced by
a $\Z_n$-grading and a nontrivial 3-cocycle.

\parindent 10pt

\section{{\Large INTRODUCTION}}

Standard methods for dealing with the nonassociativity of the
octonions involve much weaker conditions that associativity
(as alternative algebras), with the resulting problem that
most usual ideas from
linear algebra do not go through for them. In \cite{AlbMa} we have
introduced a solution to this problem based on modern ideas
from category theory and quantum group theory \cite{Ma}.
In this approach we work with algebras which {\em are} associative
but only up to a certain `rebracketting isomorphism'. A powerful
result from category theory\cite{Mac} then says that one may make all
categorical constructions exactly as for associative algebras and
afterwards insert the brackets in a consistent manner. We call such
objects {\em quasialgebras}. The paper \cite{AlbMa} particularly
studied examples of quasialgebras $k_FG$ which are obtained by
twisting from group algebras $kG$, a class that we show includes
the octonions and higher Cayley algebras. The rebracketting
isomorphism is controlled by a
3-cocycle which in these cases is a coboundary $\phi=\partial F$.
However, the theory is much more general than this and can include
other much more novel nonassociative objects. In this note we
provide examples with $\phi$ not a coboundary and hence definitely
going beyond the form $k_FG$. Our examples are related to
3-cocycles on the group $\Z_n$ and our main results include a
complete classification of the possibilities for these for low $n$.

\section{QUASIALGEBRAS}

We recall briefly the general categorical setting which we
use\cite{Ma}\cite{Ma2}. A
monoidal category $\CC$ means objects $V,W,Z$, etc equipped with a
functor $\tens:\CC\times \CC\to \CC$, and a collection of
functorial isomorphisms
\[ \Phi_{V,W,Z}:(V\tens W)\tens Z\to V\tens (W\tens Z)\]
called the rebracketting {\em associator} between any three
objects. It is required to obey the {\em pentagon identity}
\[
\begin{picture}(180,90)(115,684)
\thicklines
\put(192,697){\vector( 1, 0){ 26}}
\put(142,729){\vector( 4,-3){ 28}}
\put(253,708){\vector( 4, 3){ 28}}
\put(253,770){\vector( 4,-3){ 28}}
\put(141,749){\vector( 4, 3){ 28}}
\put(170,772){\makebox(0,0)[lb]{\raisebox{0pt}[0pt][0pt]{
$ (V\tens W)\tens (Z\tens U)$}}}
\put(222,694){\makebox(0,0)[lb]{\raisebox{0pt}[0pt][0pt]{
$ V\tens ((W\tens Z)\tens U)$}}}
\put(109,693){\makebox(0,0)[lb]{\raisebox{0pt}[0pt][0pt]{
$ ( V\tens (W\tens Z))\tens U$}}}
\put(139,761){\makebox(0,0)[lb]{\raisebox{0pt}[0pt][0pt]{
\small $\Phi$}}}
\put(285,763){\makebox(0,0)[lb]{\raisebox{0pt}[0pt][0pt]{
\small $\Phi$}}}
\put(113,713){\makebox(0,0)[lb]{\raisebox{0pt}[0pt][0pt]{
\small $ \Phi\tens\id$}}}
\put(281,711){\makebox(0,0)[lb]{\raisebox{0pt}[0pt][0pt]{
\small $ \id\tens\Phi$}}}
\put(203,684){\makebox(0,0)[lb]{\raisebox{0pt}[0pt][0pt]{
\small $\Phi$}}}
\put(100,735){\makebox(0,0)[lb]{\raisebox{0pt}[0pt][0pt]{
$ ((V\tens W)\tens Z)\tens U$}}}
\put(232,735){\makebox(0,0)[lb]{\raisebox{0pt}[0pt][0pt]{
$ V\tens (W\tens (Z\tens U))$}}}
\end{picture}
\]
which says that the two ways to reverse the brackettings as shown
coincide. Mac Lane's coherence theorem then says that all other
routes between two bracketted tensor products also coincide. In
effect, this means that one may generalise constructions in linear
algebra exactly as if $\tens$ were strictly associative, dropping
brackets. Afterwards one may add brackets, for example putting all
brackets accumulating to the left, and then insert applications of
$\Phi$ as needed for the desired compositions to make sense; all
different ways to do this will yield the same net result.

So working in such a category is no harder than usual associative
linear algebra. For example, an algebra $A$ in such a category
means
\[\bullet\circ(\bullet\tens\id)=\bullet\circ(\id\tens\bullet)
\circ\Phi_{A,A,A}\]
for the product $\bullet$, where $\Phi$ is inserted for the
bracketting to make sense. So, recognising the octonions as such a
{\em quasiassociative algebra} (or {\em quasialgebra} for short)
makes them as good as associative in the precise sense explained
above.

\begin{lemma} Let $G$ be a group
and $\phi:G\times G\times G\to k$ invertible and a cocycle:
\[ \phi(y,z,w)\phi(x,yz,w)\phi(x,y,z)=\phi(x,y,zw)\phi(xy,z,w),
\quad \phi(x,e,y)=1\]
$x,y,z,w\in G$. Then the category of $G$-graded vector spaces is
monoidal with
\[ \Phi_{V,W,Z}((v\tens w)\tens z)=v\tens
(w\tens z)\phi(|v|,|w|,|z|)\] on elements of degree
$|v|,|w|,|z|\in G$.
\end{lemma}

An algebra in this category is called a {\em $G$-graded
quasialgebra} and is by definition a $G$-graded vector space with
product respecting the grading and obeying
\[ (a\cdot b)\cdot c =a\cdot (b\cdot c)\phi(|a|,|b|,|c|),
\quad\forall a,b,c\in A\]
of homogeneous degree. There is also a notion of quasicommutativity
\[ a\cdot b=b\cdot a \CR(|a|,|b|)\]
where a quasibicharacter $\CR$ with respect to $\phi$ defines a
braiding or `generalised transposition' in the category. The
octonions are both quasiassociative and quasicommutative in the
category of $\Z_2\times\Z_2\times\Z_2$-graded spaces with
\[ \phi(\vec{x},\vec{y},\vec{z})=(-1)^{(\vec{x}\times\vec{y})
\cdot\vec{z}},\quad \CR(\vec{x},\vec{y})
=\cases{1&if\ $\vec{x}=0$\
or\ $\vec{y}=0$\ or\ $\vec{x}=\vec{y}$\cr
-1&else}\]
where we use a vector notation for the grading. Explicitly the
octonion product in the graded basis is\cite{AlbMa}
\[ e_{\vec{x}}\cdot e_{\vec{y}}
=e_{\vec{x}+\vec{y}}(-1)^{\sum_{i\le
j}x_iy_j+y_1x_2x_3+x_1y_2x_3+x_1x_2y_3}.\]

\section{ $\Z_n$-GRADED CASE}

In this section we classify the possible cocycles on $G=\Z_n$ for
low $n$ and give examples of quasialgebras of this type. We use an
additive notation throughout.

\begin{lemma}
Let $\phi$ be a cocycle on $\Z_n$ with $n$ even. The element
$x={n\over 2}$ obeys $\phi(x,x,x)=\omega$ with $\omega^2=1$.
Moreover, $\phi(x,x,y)=\omega
\phi(x,x,x+y)$ for all $y$.
\end{lemma}
\proof
Using the cocycle condition and $\phi$ trivial when any element is
the group identity, we have $\phi(x,x,x)^2={\phi(2x,x,x)
\phi(2x,x,x)\over \phi(x,2x,x)}=1$. The other result is also
immediate.
\endproof

\begin{corol} A $\Z_2$-graded quasialgebra is either an
associative superalgebra or quasiassociative with
$\phi(x,y,z)=(-1)^{xyz},
\forall x,y,z,\in \Z_2$. The latter is not a coboundary.
\end{corol}
\proof
We have $\phi(x,0,y) =\phi(0,x,y)= \phi(x,y,0)=1,\forall x,y \in
\Z_2$. By the last lemma for $x=1$ we have only two choices:
$\phi(1,1,1)=1$ or $\phi(1,1,1)=  -1$. The other result is immediate
from the formula for a coboundary. \endproof

\begin{lemma}
Let $\phi$ be a cocycle defined in $\Z_n$. Then for all $x\in \Z_n$
we have,

1) $\phi((n-1)x,x, (n-1)x)\phi(x, (n-1)x,x)=1$.

2) $\phi((n-1)x,x, x).\phi(x, x,(n-1)x)={1\over \phi((n-1)x,2x,
(n-1)x)}$.
\end{lemma}
\proof
Follows from the definition of a cocycle and $nx=0$ for all $x$.
\endproof

\begin{lemma}
Let $\phi$ be a cocycle defined in $\Z_3$.Then $\forall x\in \Z_3$
we have,

1) $\phi (2x,x,2x) \phi (x,2x,x)=1$

2) $ \phi (2x,x,x)\phi (x,x,2x)={1\over \phi (2x,2x,2x)} $

3) $(\phi (x,x,x) \phi (2x,2x,2x))^3=1$

4) $\phi (2x,2x,x) \phi (2x,x,2x)=\phi (x,x,2x) $

5) $\phi (x,x,x) \phi (x,x,2x)={\phi (2x,x,2x )\over \phi
(x,2x,2x)} $.
\end{lemma}
\proof Parts
1) and 2) follow by Lemma~4. Part 3) is $
\phi(x,x,x)^2={\phi(2x,x,x) \phi(x,x,2x)\over \phi(x,2x,x)}
$ but $  \phi(2x,x,x) \phi(x,x,2x) ={1\over \phi (2x,2x,2x)} $.
Then $  \phi(x,x,x) ^2.\phi(2x,2x,2x) \phi(x,2x,x)=1$ and
analogously
 $  \phi(2x,2x,2x) ^2.\phi(x,x,x) \phi(2x,x,2x)=1. $
So by 1) we have that $\phi(x,x,x) ^3 \phi(2x,2x,2x) ^3=1 $.
Parts 4) and 5) follow directly by the definition of a cocycle.
\endproof

\begin{propos} Every cocycle on $\Z_3$ has the form
\[ \phi_{111}=\alpha,\quad \phi_{112}=\beta,\quad \phi_{121}
={1\over \omega \alpha},\quad
\phi_{122}={\omega \over \beta}\]
\[\phi_{211}={\alpha\over\beta\omega},\quad \phi_{212}
={\alpha \omega},\quad
\phi_{221}={\beta \over\omega\alpha}, \quad
\phi_{222}= {\omega \over \alpha}\]
for some non zero $\alpha,\beta \in k$ and $\omega$ a cubic
root of the unity. Here $\phi(1,1,1)=\phi_{111}$, etc. is a
shorthand.
\end{propos}
\proof First of all, let
$\omega=\phi(1,1,1)\phi(2,2,2)$, a cubic root of unity by part 3)
of the last lemma. We also have $\phi(1,1,1) ^2={
\phi(2,1,1)
\phi(1,1,2)
\over
\phi(1,2,1)}={1\over  \phi(1,2,1)\phi(2,2,2) } $ by part 2) of
the last lemma. Hence $\phi(1,1,1)={1\over \omega \phi(1,2,1)} $.
On the other hand $\phi(1,1,2)\phi(1,2,2)=\phi(2,2,2)
\phi(1,1,1)=\omega$. Also from part 2) of the lemma, we have
$\phi(2,1,1)\phi(1,1,2)={1\over\phi(2,2,2)}$ and from part 4) we
have $\phi(2,2,1)\phi(2,1,2)=
\phi(1,1,2)$. Similarly part 5) gives $\phi(1,1,1)\phi(1,1,2)
={\phi(2,1,2)\over
\phi(1,2,2)}$. Denoting
$\phi(1,1,1)=\alpha$ and $\phi(1,1,2)=\beta$, we have the result as
stated. \endproof

This can be written, for example, as
\[ \phi(x,y,z)=(\alpha^{(-1)^z+x-xz}\beta^{x-z})^{(-1)^y}
\cases{1&{\rm
if}\ $x=y=1$\cr \omega^z&{\rm else}}\]
for $x,y,z\ne 0$.

\begin{propos} A cocycle on $\Z_3$ in the parametrisation
above is coboundary iff  $\omega=1$.
\end{propos}
\proof Chose any cochain $F$ (an invertible function such that
$F(0,x)=F(x,0)=1$
for all $x$). For brevity we write its entries as a matrix
$F(1,2)=F_{12}$ etc. Let ${ F_{21}\over F_{12} }=\alpha$ and ${
F_{11}F_{22}\over F_{12}}=\beta$. Then a coboundary
$\phi(x,y,z)=F(x,y)F(xy,z)/F(y,z)F(x,yz)$ is  $\phi_{111}
={ F_{11}F_{21}\over F_{11}F_{12} }=\alpha$,
$\phi_{112}={ F_{11}F_{22} \over F_{12}}=\beta$, $\phi_{121}={
F_{12} \over F_{21} }={1\over \alpha}$, $\phi_{122}={ F_{12}\over
F_{22}F_{11}}={1\over\beta}$, $\phi_{211}={ F_{21}\over
F_{11}F_{22} }={\alpha\over\beta}$, $\phi_{212}={ F_{21}\over
F_{12}}=\alpha,\quad \phi_{221}={ F_{22}F_{11}\over
F_{21}}={\beta\over\alpha}$, $\phi_{222}={ F_{22}F_{12} \over
F_{22}F_{21} }={1\over\alpha}$ which is of the form above with
$\omega=1$. Conversely, if $\phi$ of the form above is a coboundary
then $\phi_{111}={ F_{11}F_{21}
\over F_{11}F_{12} }=\alpha,
\phi_{121}={ F_{12}F_{01}\over F_{21}F_{10}}
={1\over \omega\alpha}$. So $\omega=1$.
\endproof

\begin{propos} Every choice of invertible $\alpha,\beta,
\omega$ with $\omega^3=1$
yields a cocycle on $\Z_3$ of the form above. In particular,
\[ \phi(x,y,z)=\cases{1&if\ $x=y=1$\cr \omega^z&else}\]
for $x,y,z\ne 0$ is a
noncoboundary cocycle when $\omega\ne 1$ and every cocycle is
cohomologically equivalent to one of this form.
\end{propos}
\proof  We take $\alpha=\beta=1$ and $\omega$ a nontrivial cube
root of unity in Proposition~6 and verify directly that it is
indeed a 3-cocycle. The cocycle condition is empty when any of the
arguments is $0$, so we assume that they are not. Then, as we
have two different expressions for $\phi$, we consider the cases
(i) $\phi(1,1,z) \phi(1,z,w)={\phi(2,z,w)
\phi(1,1,z +w)\over\phi(1,1 +z,w)}$
is satisfied because both sides are $1$ if $z=1$ and
$\omega^w$ when $z=2$.
(ii) $\phi(x,1,1) \phi(1,1,w)={\phi(x+ 1,1,w)
\phi(x,1,1 +w)\over\phi(x,2,w)}$
is satisfied because both sides are $1$ if $x=1$ and $\omega$
if $x=2$. (iii) $\phi(x,1-x,1) \phi(1-x,1,w)={\phi(1,1,w)
\phi(x,1-x,1+w)\over\phi(x,2-x,w)}$ is satisfied because both sides
are $1$ if $x=1$ and $\omega^{w+1}$ if $x=2$. (iv)
$\phi(1,y,1-y) \phi(y,1-y,w)={ \phi(1+y,1-y,w) \phi(1,y,1-y
+w)\over\phi(1,1,w)}$ is satisfied because
both sides are $1$ if $y=1$ and $\omega^{w-1}$ if $y=2$.
On the other hand, we know by Proposition~6 that every cocycle
is the product of this one defined by some $\omega$
and one of the coboundary type defined by $\alpha,\beta$ in
Proposition~7.
\endproof

A more symmetric choice to generate the cohomology is
with $\alpha=\beta=\omega^2$. This can be written more
compactly as
\[ \phi(x,y,z)=\omega^{xz-xy-yz}\]
for $x,y,z\ne 0$, and is cohomologically equivalent to
the cocycle in
Proposition~8.

\begin{corol} A cocycle on $\Z_3$ is trivial if and
only if there is an element $x\ne 0$ in $\Z_3$ such that
$\phi(x,y,z)=1$ for all $y,z$
\end{corol}
\proof  If $\phi (1,y,z)=1$ for all $y,z$, we have
$\phi_{112}=\phi_{111}=\phi_{121}=1$ and hence by
Proposition~6 we have
$\alpha=\beta=\omega=1$ and $\phi=1$.
If $\phi (2,y,z)=1$ for all $y,z$ we have that
$\phi_{212}=\phi_{211}=\phi_{222}=1$ and hence
$\alpha=\beta=\omega=1$ again.
\endproof

Returning to the general case, a natural cocycle motivated by
some of the above is:

\begin{corol} Let $q$ be an $n$-th root of unity. Then
\[ \phi(x,y,z)=q^{xyz}\]
is a cocycle on $\Z_n$. When $n=3$ it is a coboundary with
$\alpha=q,\beta=q^2,\omega=1$.
\end{corol}
\proof The 3-cocycle condition becomes $yzw+x(y+z)w+xyz
=xy(z+w)+(x+y)zw$
in $\Z_n$ and holds using distributivity of the product
in the ring $\Z_n$ over the additive group structure.
When $n=3$ it must fit into our classification above,
which it does with $\omega=1,\alpha=q,\beta=q^2$. \endproof

Note that if $q$ is not a root of unity, we still have a
coboundary,
\[ \phi_{111}=\phi_{212}=\phi_{221}=q,\quad \phi_{121}
=\phi_{211}=\phi_{222}=q^{-1},\quad
\phi_{122}=q^{-2},\quad\phi_{112}=q^2\]
if we choose $\alpha=q,\beta=q^2$ again.

Let us stress that every cocycle leads to a category of
quasialgebras and that these are different even if the
cocycles are cohomologically equivalent, i.e we are interested
in the full parametrisation in Proposition~6.
When related by a coboundary the quasialgebras may potentially
be related to each other by twisting in the same way as the
octonions are a twist of the group algebra of $\Z_2\times
\Z_2\times \Z_2$\cite{AlbMa}.
When in different cohomology classes then the quasialgebras cannot
 be related by a twist and are in this sense
`topologically distinct' examples.

\begin{example} Let $q$ be a cubic root of unity. The $\Z_3$-graded
quasialgebra with graded basis $\{e_x\}$
for $x\in\Z_3$, $e_0=1$ and other products
\[ e_x e_y=e_{x+y}q^{y-x},\quad\forall x,y\ne 0\]
is a twisting of $k\Z_3$ and has coboundary cocycle
$\phi(x,y,z)=q^{xyz}$.
\end{example}
\proof We take $F(x,y)=q^{y-x}$ for $x,y\ne 0$ in Proposition~7.
This has $\alpha=q,\beta=q^2$ and hence $\partial F=\phi$ in
Corollary~10. \endproof

On the other hand, Theorem~7.3 of \cite{AlbMa} provides
 a construction of a quasialgebra for {\em any} cocycle (and any
 graded vector space) as the quasialgebra of quasi-matrices. In
particular, let $\phi$ be a cocycle on $\Z_n$ then the natural
quasialgebra of quasimatrices $M_{n,\phi}$ has basis
$E_{ij}$ labelled by $i,j\in\Z_n$ and of degree $i-j$, with
the product
\eqn{qmat}{ E_{ij}\cdot E_{kl}=\delta_{jk}E_{il}
{\phi(i,-j,j-l)\over\phi(-j,j,-l)}.}
The quasiassociativity is
\eqn{qmatassoc}{ (E_{ij}\cdot E_{kl})\cdot E_{rs}
=E_{ij}\cdot (E_{kl}\cdot E_{rs})\phi(i-j,k-l,r-s)}
which can be computed more explicitly depending on the
form of $\phi$.

\begin{example} Let $q^n=1$. Then $M_{n,\phi}$
for $\phi$ in Corollary~10 has the product
\[ E_{ij}\cdot E_{kl}=\delta_{jk}E_{il}q^{ijl-j^2(i+l)}.\]
Let $\omega^3=1$. Then $M_{3,\phi}$ for the noncoboundary
$\phi$ in Proposition~8 has the product
\[ E_{ij}\cdot E_{kl}
=\delta_{jk}E_{il}
\cases{\omega^l& {\rm if}\ $i=0$, $j\ne0$\ {\rm or}\ $i=1$, $j=2$\cr
\omega^j&{\rm else}.}\]
\end{example}
\proof We insert the form of the relevant cocycle into
(\ref{qmat}). In
the second case all the possibilities for $i,j,l,j-l$
zero or not have to be looked at separately but can
afterwards be recombined as stated. \endproof

Here $M_{2,\phi}$ for $q\ne 1$ is also noncoboundary and
an example of the second type in Corollary~3. Also, the
corresponding quasimatrix product\cite{AlbMa} among actual
matrices $a,b$ has the same form
\eqn{qmatmul}{ (a\cdot b)_{il}
=\sum_j a_{ij}b_{jl} {\phi(i,-j,j-l)\over\phi(-j,j,-l)}}
and therefore the same coefficients for the $M_{n,\phi}$,
$M_{3,\phi}$ as appearing in Example~12. For example,
$M_{2,\phi}$ has the product
\[ \pmatrix{a_{00}&a_{01}\cr a_{10}&a_{11}}
\cdot\pmatrix{b_{00}&b_{01}\cr b_{10}&b_{11}}
=\pmatrix{a_{00}b_{00}+a_{01}b_{10}&
a_{00}b_{01}-a_{01}b_{11}\cr
a_{10}b_{00}-a_{11}b_{10}&a_{10}b_{01}-a_{11}b_{11}}.\]

Finally, we note that the cocycle in Corollary~10 has an obvious
generalisation to $(\Z_n)^m$ as
\[ \phi(\vec{x},\vec{y},\vec{z})=q^{(\vec{x},\vec{y},\vec{z})}\]
where we use a vector notation with components in $\Z_n$ and $(\ ,\
,\ )$ is $\Z_n$-trilinear. The cocycle for the octonions is a
coboundary example of this type on $(\Z_2)^3$.

\end{document}